

\baselineskip=14pt
\parskip=10pt
\def\halmos{\hbox{\vrule height0.15cm width0.01cm\vbox{\hrule height
  0.01cm width0.2cm \vskip0.15cm \hrule height 0.01cm width0.2cm}\vrule
  height0.15cm width 0.01cm}}
\font\eightrm=cmr8 

\magnification=\magstephalf

\def\1{{\overline{1}}}
\def\2{{\overline{2}}}
\parindent=0pt
\overfullrule=0in

\def\frac#1#2{{#1 \over #2}}

\centerline
{\bf The Dayenu Boolean Function Is Almost Always True! }
\bigskip
\centerline
{\it Doron ZEILBERGER}

{\qquad} {\it Dedicated to Hemi and Yael Nae}

[{\eightrm Recall, from Logic 101,  that for any two statements, $X,Y$, $X \vee Y$  means that ``either $X$ or $Y$ or both happened'',
$X \wedge Y$ means that {\bf both} $X$ and $Y$ happened, while $\bar{X}$ means that $X$ did {\bf not} happen.
For typographical clarity, $X \wedge Y$ is often written $XY$.}]

Two nights ago, my wife Jane and I were fortunate to be guests in a wonderful Passover {\it seder} at the house
of Hemi\footnote{$^1$}{
{
\eightrm Hemi and I were dorm-mates at the Weizmann Institute, way back in the early seventies. His son is also called Doron,
so we call him `little Doron', while I am `big Doron'. `Little Doron' is no longer so little (he is $25$-years-old),
but Hemi commented that he knows me (`big Doron') much longer than he knows `little Doron'.}
} and Yael Nae. Soon enough we came to the number-one hit song, {\it Dayenu},
praising God for doing $15$ amazing miracles, let's call them $x_1, \dots, x_{15}$, where
$x_1$ stands for ``Took us out of Egypt'', and $x_{15}$ stands for  ``Built us our Temple''.
The exact nature of the other $13$ miracles (and for that matter the above two ) are irrelevant
to this {\it mathematical} paper, but can be easily looked-up in any {\it Haggadah} (and of course, nowadays, on the internet).

It so happened that what God actually did can be described by the Boolean function (with $n=15$)
$$
G_n(x_1, \dots, x_n)= x_1 \wedge \dots  \wedge x_n = \bigwedge_{i=1}^{n} x_i \quad,
$$
whose only {\it truth} vector is the all-true vector $T^n$, i.e., assuming, that the probability of
any one miracle occurring is $p$ (and that they are independent events), 
has the tiny probability of $p^n$.

But the author of {\it Dayenu} asserts that God was an over-achiever, and we, the children of Israel, should have
been content if the following Boolean function would have been {\it satisfied}.

$$
D_n(x_1, \dots, x_n) = \bigvee_{i=1}^{n-1} x_i \,\, \overline{x_{i+1}}  \quad  .
$$

The following question {\it immediately} came to my mind: How many (and which) truth-vectors are satisfied by
the {\it Dayenu} function, and what is the probability that God, deciding randomly which miracles to perform and
which not to perform, would have satisfied the {\it minimum} requirement demanded by the anonymous author of {\it Dayenu}?

Thanks to De Morgan we have
$$
\overline{D_n(x_1, \dots, x_n)} = \bigwedge_{i=1}^{n-1} \,\, (\, \overline{x_i} \, \vee  \, x_{i+1} \, ) \quad .
$$

We are now ready for 

{\bf The Dayenu Theorem.} The {\it} full disjunctive normal form of the negation of the Dayenu Boolean function
is given by
$$
\overline{D_n(x_1, \dots, x_n)} = 
\bigvee_{i=1}^{n+1} \quad \left ( \, \bigwedge_{j=1}^{n-i+1} \overline{x_j}\bigwedge_{j=n-i+2}^{n} x_j  \, \right )\quad ,
$$
or, more concretely:
$$
\overline{D_n(x_1, \dots, x_n)} = 
\overline{x_1} \, \overline{x_2} \, \cdots  \overline{x_{n-1}}  \, \overline{x_{n}} \,\, \vee \,\,
\overline{x_1} \, \overline{x_2} \, \cdots  \overline{x_{n-1}}\,x_n \,\,  \vee \,\,
\, \dots  \, \, \vee \,\,
\overline{x_1} \, x_2 \, \dots   x_{n-1}\,x_n  \, \,  \vee \, \, 
x_1 \, x_2 \, \cdots \, x_n \quad .
$$

{\bf Proof}: By induction on $n$. It is true when $n=2$ (check!). 
Assume that it is true when $n$ is replaced by $n-1$.
Note that
$$
\overline{D_n(x_1, \dots, x_n)} = 
 \overline{D_{n-1}(x_1, \dots, x_{n-1})} \,\wedge \, \,  \left ( \overline{x_{n-1}} \, \vee \, x_n \right ) 
$$
$$
=
 \overline{D_{n-1}(x_1, \dots, x_{n-1})} \,\wedge \, \overline{x_{n-1}} \, \, \vee \, \, \overline{D_{n-1}(x_1, \dots, x_{n-1})} \, \wedge \, x_n  \quad .
$$

By the  inductive hypothesis:
$$
\overline{D_{n-1}(x_1, \dots, x_{n-1})} \, =  \,
\overline{x_1} \, \overline{x_2} \, \cdots  \overline{x_{n-2}}  \, \overline{x_{n-1}} \, \, \vee \, \,
\overline{x_1} \, \overline{x_2} \, \cdots  \overline{x_{n-2}}\,x_{n-1}  \,\, \vee \,\,
\, \dots  \, \, \vee \, \,
\overline{x_1} \, x_2 \, \dots   x_{n-2}\,x_{n-1} \,\,  \vee \, \,
x_1 \, x_2 \, \cdots \, x_{n-1} \quad .
$$
Regarding  $\overline{D_{n-1}(x_1, \dots, x_{n-1})} \, \wedge  \, \overline{x_{n-1}}$ we have
$$
 \overline{D_{n-1}(x_1, \dots, x_{n-1})} \, \wedge \, \overline{x_{n-1}} 
$$
$$
=
\overline{x_1} \, \overline{x_2} \, \cdots  \overline{x_{n-2}}  \, \overline{x_{n-1}} \,  \overline{x_{n-1}} \, \, \vee \, \,
\overline{x_1} \, \overline{x_2} \, \cdots  \overline{x_{n-2}}\,x_{n-1}  \, \overline{x_{n-1}}\, \,  \vee \,\, \dots  \, \, \vee \,\,
\overline{x_1} \, x_2 \, \dots   x_{n-2}\,x_{n-1}  \, \overline{x_{n-1}} \vee \,
x_1 \, x_2 \, \cdots \, \, x_{n-1} \, \overline{x_{n-1}} \quad .
$$
Since  $x_{n-1}  \, \overline{x_{n-1}}$ is FALSE, all the above terms, except the first, vanish, 
and since  $\overline{x_{n-1}}  \, \overline{x_{n-1}}= \overline{x_{n-1}}$, we have
$$
 \overline{D_{n-1}(x_1, \dots, x_{n-1})} \, \wedge \, \overline{x_{n-1}} =
\overline{x_1} \, \overline{x_2} \, \cdots  \overline{x_{n-2}}  \, \overline{x_{n-1}} \quad .
$$
Regarding  $\overline{D_{n-1}(x_1, \dots, x_{n-1})} \,\wedge \, x_n$ we have
$$
\overline{D_{n-1}(x_1, \dots, x_{n-1})} \, x_n  \, =
$$
$$
\overline{x_1} \, \overline{x_2} \, \cdots  \overline{x_{n-2}}  \, \overline{x_{n-1}} \, x_n \, \, \vee \, \,
\overline{x_1} \, \overline{x_2} \, \cdots  \overline{x_{n-2}}\,x_{n-1} \, x_n \,\,  \vee \,\,
 \dots  \, \, \vee \,\,
\overline{x_1} \, x_2 \, \dots   x_{n-2}\,x_{n-1}  \,  x_n \,\, \vee \, \,
x_1 \, x_2 \, \cdots \, x_{n-1}  \, x_n \quad .
$$
Combining, we have
$$
\overline{D_n(x_1, \dots, x_n)}  \, =
$$
$$
\overline{x_1} \, \overline{x_2} \, \cdots  \overline{x_{n-2}}  \, \overline{x_{n-1}} 
$$
$$
\, \, \vee \, \,
\overline{x_1} \, \overline{x_2} \, \cdots  \overline{x_{n-2}}  \, \overline{x_{n-1}} \, x_n \, \, \vee \, \,
\overline{x_1} \, \overline{x_2} \, \cdots  \overline{x_{n-2}}\,x_{n-1} \, x_n  \, \, \vee \,\,
\, \dots  \, \vee \,\,
\overline{x_1} \, x_2 \, \dots   x_{n-2}\,x_{n-1}  \,  x_n \,\, \vee \, \,
x_1 \, x_2 \, \cdots \, x_{n-1}  \, x_n \quad .
$$
Since
$$
\overline{x_1} \, \overline{x_2} \, \cdots  \overline{x_{n-2}}  \, \overline{x_{n-1}} \, = \,
\overline{x_1} \, \overline{x_2} \, \cdots  \overline{x_{n-2}}  \, \overline{x_{n-1}} \, \overline{x_{n}} \, \, \vee \,\,
\overline{x_1} \, \overline{x_2} \, \cdots  \overline{x_{n-2}}  \, \overline{x_{n-1}} \, x_{n}  \quad ,
$$
we get
$$
\overline{D_n(x_1, \dots, x_n)}  \, = 
$$
$$
\overline{x_1} \, \overline{x_2} \, \cdots  \overline{x_{n-2}}  \, \overline{x_{n-1}} \, \overline{x_{n}} \,\, \vee \,\,
\overline{x_1} \, \overline{x_2} \, \cdots  \overline{x_{n-2}}  \, \overline{x_{n-1}} \, x_{n} \, \, \vee 
$$
$$
\overline{x_1} \, \overline{x_2} \, \cdots  \overline{x_{n-2}}  \, \overline{x_{n-1}} \, x_n \, \,\vee \,\,
\overline{x_1} \, \overline{x_2} \, \cdots  \overline{x_{n-2}}\,x_{n-1} \, x_n  \, \, \vee \,\,
\, \dots  \, \, \vee \, \,
\overline{x_1} \, x_2 \, \dots   x_{n-2}\,x_{n-1}  \,  x_n \,\, \vee \, \,
x_1 \, x_2 \, \cdots \, x_{n-1}  \, x_n \quad .
$$
Since the second and third term above are the same (and $X \vee X=X$), we finally get
$$
\overline{D_n(x_1, \dots, x_n)}  \, =
$$
$$
\overline{x_1} \, \overline{x_2} \, \cdots  \overline{x_{n-2}}  \, \overline{x_{n-1}} \, \overline{x_{n}} \,\, \vee \,\,
\overline{x_1} \, \overline{x_2} \, \cdots  \overline{x_{n-2}}  \, \overline{x_{n-1}} \, x_n \,\, \vee \,\,
\overline{x_1} \, \overline{x_2} \, \cdots  \overline{x_{n-2}}\,x_{n-1} \, x_n  \, \, \vee \, \,
\dots \, \, \vee \,\,
\overline{x_1} \, x_2 \cdot x_n  \, \, \vee \,\,
x_1 \, x_2 \, \cdots \, x_{n-1}  \, x_n \quad .
$$
\halmos

{\bf Corollary}: Assuming that the probability of each miracle is $\frac{1}{2}$, and that they are done
independently, the probability of {\bf not} meeting the Dayenu requirement is $\frac{n+1}{2^n}$ and
hence of meeting it is $1-\frac{n+1}{2^n}$, that happens to be, for $n=15$,
$\frac{2047}{2048}=  0.9995117\dots$ \quad .

{\bf Comments}:

{\bf 1.} Surprisingly, what God actually did, performing {\it all} the miracles, is not part of the
truth-set of the Dayenu Boolean function, since  there is always at least one miracle that is {\bf not} performed.

{\bf 2.} A faster proof, without Boolean logic, for getting the set of true-false vectors not satisfying the
Dayenu function $D_n$ (and hence proving the above Dayenu theorem),
can be gotten by finding the set of all members of $\{T,F\}^n$
for which an $F$ {\bf never} (immediately) follows a $T$. 
Of course, this set is $\{ F^n,F^{n-1}T, \dots, T^n\}$.

\bigskip
\hrule
\bigskip

{\bf 16 Nisan 5777}. Exclusively published in the Personal Journal of Shalosh B. Ekhad and Doron Zeilberger ,
{\tt http://www.math.rutgers.edu/\~{}zeilberg/pj.html}, and arxiv.org .

\bigskip
\hrule
\bigskip
Doron Zeilberger, Department of Mathematics, Rutgers University (New Brunswick), Hill Center-Busch Campus, 110 Frelinghuysen
Rd., Piscataway, NJ 08854-8019, USA. \hfill\break
Email: {\tt DoronZeil at gmail  dot com}   \quad .
\bigskip
\hrule
\bigskip

\end